\documentclass[12pt]{article}
\input isolatin1.sty
 \usepackage{a4}
 \usepackage{amsmath}
 \usepackage{amssymb}
\usepackage{amsthm}
\usepackage{graphicx}
\newtheorem{thm}{Theorem}[section]

\newtheorem{lem}[thm]{Lemma}

\newtheorem{rem}[thm]{Remark}

\numberwithin{equation}{section} %%%%%%numera las ecuaciones (1.1), en
\title{Iteration of the rational function $z-1/z$ and a Hausdorff
  moment sequence
\footnote{This work was initiated during a visit of the first author
  to the University
of Sevilla in January 2007 supported by  D.G.E.S,
  ref. BFM2003-6335-C03-01, FQM-262 ({\em Junta de Andaluc\'{\i}a}).}}
\author{Christian Berg $^{\dagger}$ and Antonio J. Dur\'{a}n  $^{\ddagger}$\\
   \footnotesize $\dagger$ \ Institut for Matematiske Fag. K\o benhavns
   Universitet \\
    \footnotesize Universitetsparken 5; DK-2100 K\o benhavn \o,
   Denmark. berg@math.ku.dk \\  \footnotesize
   $\ddagger$ \  Departamento de An\'{a}lisis Matem\'{a}tico.
   Universidad de Sevilla \\
   \footnotesize Apdo (P. O. BOX) 1160. 41080 Sevilla. Spain. duran@us.es \\
\ \ }

 \date{\today}
\begin{document}
\maketitle

\begin{abstract} In a previous paper we considered a positive function $f$,
  uniquely determined for $s>0$ by the requirements $f(1)=1$,
  $\log(1/f)$ is convex and the functional equation
  $f(s)=\psi(f(s+1))$ with $\psi(s)=s-1/s$. We prove that the
  meromorphic extension of $f$ to the whole complex plane is given by
  the formula $f(z)=\lim_{n\to\infty}\psi^{\circ
    n}(\lambda_n(\lambda_{n+1}/\lambda_n)^z)$, where the numbers
  $\lambda_n$ are defined by $\lambda_0=0$ and the recursion $\lambda_{n+1}=(1/2)(\lambda_n+\sqrt{\lambda_n^2+4})$. 
The numbers $m_n=1/\lambda_{n+1}$ form a Hausdorff moment sequence of
a probability measure $\mu$ such that $\int t^{z-1}\,d\mu(t)=1/f(z)$.
\end{abstract}



\section{Introduction and main results}
 Hausdorff moment sequences are sequences
of the form $a_n=\int _0^1t^nd \nu (t)$, $n\ge 0$, where $\nu $ is a
positive measure on $[0,1]$. 

  In \cite{bd2} we introduced a non-linear transformation
$T$ of the set of Hausdorff moment sequences into itself
by the formula:
\begin{equation}\label{eq:tms}
T((a_n))_n=1/(a_0+a_1+\cdots +a_n),\quad n\ge 0.
\end{equation}
The corresponding transformation of  positive measures on $[0,1]$ is 
denoted $\widehat{T}$. We recall from \cite{bd2} that if $\nu\neq 0$, then
$\widehat{T}(\nu)(\{0\})=0$ and
\begin{equation}\label{eq:bd2}
\int_0^1\frac{1-t^{z+1}}{1-t}\,d\nu(t)\int_0^1 t^z\,d\widehat{T}(\nu)(t)=1
\;\text{for}\; \Re z \geq 0.
\end{equation}

It is clear that if $\nu $ is a probability measure, then
so is $\widehat{T}(\nu )$, and in this way we get a transformation of the
convex set of normalized Hausdorff moment sequences (i.e. $a_0=1$) as well as a
transformation of the set of probability measures on $[0,1]$. 
By Kakutani's theorem the transformation has a fixed point, and by
(\ref{eq:tms}) it is clear that a fixed point $(m_n)_n$ is uniquely
 determined by the
recursive equation
\begin{equation}\label{eq:fix}
m_0=1,\quad (1+m_1+\cdots +m_n)m_n=1,\quad n\ge 1.
\end{equation}                                    
Therefore
\begin{equation}\label{eq:fixrec}
m_{n+1}^2+\frac{m_{n+1}}{m_n}-1=0,
\end{equation}
giving
$$
m_1=\frac{-1+\sqrt 5}{2},\quad m_2=\frac{\sqrt{22+2\sqrt 5}-\sqrt
5-1}{4}, \cdots\,.
$$

In \cite{bd3} we studied the Hausdorff moment sequence
$(m_n)_n$
and its associated probability measure $\mu $, called the
{\it fixed point measure}. It has an increasing and
convex density $\mathcal D$ with respect to Lebesgue measure on
$]0,1[$, and for $t\to 1$ we have $\mathcal D(t)\sim
1/\sqrt{2\pi(1-t)}$.

 We  studied the
 {\it Bernstein transform}
\begin{equation}\label{eq:Bern}
f(z)=\mathcal B(\mu)(z)=\int _0^1\frac{1-t^z}{1-t}d\mu (t),\quad  \Re z> 0
\end{equation}
as well as the {\it Mellin transform}
\begin{equation}\label{eq:Mellin}
F(z)=\mathcal M(\mu)(z)=\int _0^1t^zd\mu (t),\quad  \Re z> 0
\end{equation}
of $\mu$.
These functions are clearly holomorphic in
the half-plane $\Re z>0$ and continuous in $\Re z\ge 0$, the latter because
$\mu(\{0\})=0$.

As a first result we proved:

\begin{thm}[\cite{bd3}]\label{thm:main-a} The functions $f,F$ can be
 extended to
 meromorphic functions in $\mathbb C$ and they satisfy
\begin{equation}\label{eq:fF=1}
f(z+1)F(z)=1,\quad z\in\mathbb C
\end{equation}
\begin{equation}\label{eq:z,z+1}
f(z)=f(z+1)-\frac{1}{f(z+1)},\quad z\in\mathbb C.
\end{equation}
They are holomorphic in $\Re z>-1$. Furthermore $z=-1$ is
a simple pole of $f$ and $F$.

The fixed point measure $\mu$ has the properties
\begin{equation}\label{eq:Landau}
\int_0^1 t^x\,d\mu(t)<\infty,\;\; x>-1;\quad 
\int_0^1 \frac{d\,\mu(t)}{t}=\infty.
\end{equation}
\end{thm}

The function $f$ can be characterized in analogy with the
Bohr-Mollerup theorem about the Gamma function, cf. \cite{ar}. More
precisely we proved:

\begin{thm}[\cite{bd3}]\label{thm:main-b} The Bernstein transform (\ref{eq:Bern})
of the fixed point measure is a function  $f:\left]0,\infty\right[\to \left]0,\infty\right[$
 with the following properties
\begin{enumerate}
\item[(i)] $f(1)=1$,
\item[(ii)] $\log(1/f)$ is convex,
\item[(iii)] $f(s)=f(s+1)-1/f(s+1),\quad s>0$.
\end{enumerate}
Conversely, if $\tilde f:]0,\infty[\to ]0,\infty[$ satisfies
(i)-(iii), then it is equal to $f$ and
for $0<s\le 1$ we have
\begin{equation}\label{eq:iterate}
\tilde f(s)=\lim_{n\to\infty}\psi^{\circ n}\left(\frac{1}{m_{n-1}}\left(
\frac{m_{n-1}}{m_n}\right)^s\right),
\end{equation}
where $\psi$ is the rational function
\begin{equation}\label{eq:psi}
\psi(z)=z-\frac{1}{z}.
\end{equation}
In particular 
(\ref{eq:iterate}) holds for $f$.
\end{thm}

In this paper we use the notation for composition of mappings:
$$
\psi^{\circ 1}(z)=\psi(z), \psi^{\circ
  n}(z)=\psi(\psi^{\circ(n-1)}(z)),\;n\ge 2.
$$

The function $\psi$ is a two-to-one mapping of $\mathbb
C\setminus\{0\}$ onto $\mathbb C$ with the exception that $\psi(z)=\pm
2i$ has only one solution $z=\pm i$. Moreover,
$\psi(0)=\psi(\infty)=\infty$ and $\psi$ is a continuous mapping of
the Riemann sphere $\mathbb C^*=\mathbb C\cup \{\infty\}$ onto itself.
 It is strictly increasing on the
half-lines $\left]-\infty,0\right[$ and $\left]0,\infty\right[$,
mapping each of them onto $\mathbb R$.
The functional equation (\ref{eq:z,z+1}) can be written
\begin{equation}\label{eq:it}
f(z)=\psi(f(z+1)).
\end{equation}
The sequence $(\lambda _n)_n$ is defined in terms of $(m_n)_n$ from
(\ref{eq:fix}) by
\begin{equation}\label{eq:sl}
\lambda _0=0, \quad \lambda _{n+1}=1/m_n,\quad n\ge 0,
\end{equation}
i.e.
$$
\lambda_1=1,\quad\lambda_2=\frac{1+\sqrt{5}}2,\quad \lambda_3=\frac{\sqrt{22+2\sqrt 5}+\sqrt
5+1}{4},\cdots.
$$
By (\ref{eq:Mellin}) and (\ref{eq:fF=1}) we clearly have 
\begin{equation}\label{eq:m-la}
m_n=F(n),\;\lambda_n=f(n),\quad n\ge 0,
\end{equation}
hence by (\ref{eq:it})
\begin{equation}\label{eq:fes1}
\lambda_{n}=\psi(\lambda_{n+1}),\quad n\ge 0,
\end{equation}
which can be reformulated to
\begin{equation}\label{eq:fes2}
 \lambda _{n+1}=\frac{1}{2}\left( \lambda _{n}+\sqrt {\lambda _n^2+4}\right),
\quad n\ge 0.
\end{equation}

The main purpose of this paper is to prove that equation
\eqref{eq:iterate} holds for $f$ in the whole complex plane.

For a domain $G\subseteq\mathbb C$ we denote by $\mathcal H(G)$ the
space of holomorphic functions on $G$ equipped with the topology of
uniform convergence on compact subsets of $G$.

\begin{thm}\label{thm:main-c} Let $a_n,b_n:\mathbb C\to\mathbb C$
  be the entire functions defined by
\begin{equation}\label{eq:anbn}
a_n(z)=\lambda _n \left(\frac{\lambda _{n}}{\lambda _{n-1}}\right)
^z,n\ge 2,\quad
b_n(z)=\lambda _n\left( \frac{\lambda _{n+1}}{\lambda _{n}}\right)
^z,n\ge 1.
\end{equation}
The meromorphic function $f$ from Theorem \ref{thm:main-a} is given
for $z\in\mathbb C$ by
\begin{equation}\label{eq:iterateab}
f(z)=\lim_{n\to\infty}\psi^{\circ n}\left(a_n(z)\right)=
\lim_{n\to\infty}\psi^{\circ n}\left(b_n(z)\right),
\end{equation}
and the convergence is uniform on compact subsets of $\mathbb C$.
\end{thm}

\begin{rem} {\rm In \cite{bd3} it was proved that the Julia and Fatou
    sets of $\psi$ are respectively $\mathbb R^*$ and $\mathbb
    C\setminus\mathbb R$. For $z\in\mathbb C\setminus\mathbb R$ we
    have $\psi^{\circ n}(z)\to\infty$ for $n\to\infty$. Notice that
    $a_n(z),b_n(z)$ are close to $\lambda_n$ when $n$ is large because
$\lambda_{n+1}/\lambda_n\to 1$ according to Lemma \ref{thm:aps}
below. Also $\psi^{\circ n}(\lambda_n)=0$ for all $n$.}
\end{rem}

Hausdorff's characterization of moment sequences was given in
\cite{ha}. See also Widder's monograph \cite{W}. For information about
moment sequences in general see \cite{ak}.
\section{Proofs}

We first recall some properties of the sequence $(\lambda
_n)_n$, which are needed in the proof of Theorem \ref{thm:main-c}.

\begin{lem}[\cite{bd3}]\label{thm:aps}
\begin{enumerate}
\item $\displaystyle \sqrt n\le \lambda _n\le \sqrt{2n}$, $n\ge 0$.
\item $(\lambda_n)_n$ is an increasing divergent sequence
and $\lambda_{n+1}/\lambda_n$ is decreasing with
$\displaystyle \lim _{n\to\infty} \frac{\lambda _{n+1}}{\lambda _n}=1$.
\item $\displaystyle \lim _{n\to\infty} (\lambda ^2_{n+1}-\lambda ^2_n)=2$.
\item $\displaystyle \lim _{n\to\infty} \frac{\lambda ^2_n}{n}=2$.
\end{enumerate}
\end{lem}

{\sl Proof of Theorem \ref{thm:main-c}}.

\medskip

The proof is given in a number of steps.

\medskip
\noindent $1^{\circ}:$ {\it For any $0<s<\infty$ we have}
$$
\lim_{n\to\infty}\left[\psi^{\circ n}(a_n(s))-\psi^{\circ n}(b_n(s))\right]=0.
$$

\medskip
Note that the quantity under the limit is positive because
$a_n(s)>b_n(s)$ and $\psi$ is increasing. 

By the mean value theorem we get for a  certain
$w\in ]b_n(s),a_n(s)[$ 
$$
\psi ^{\circ n} (a_n(s))-\psi ^{\circ n} (b_n(s))=
(a_n(s)-b_n(s))(\psi ^{\circ n})' (w)
$$
$$
=(a_n(s)-b_n(s))\psi '(\psi ^{\circ n-1} (w))\psi '(\psi ^{\circ n-2} (w))
\cdots \psi '(w).
$$

Since $\lambda _n<b_n(s)<w<a_n(s)$, we get
$\lambda _{n-k}<\psi ^{\circ k}(b_n(s))<\psi ^{\circ k}(w)$, $k=0,1,
\ldots ,n$, hence
\begin{eqnarray*}
\lefteqn{\vert \psi ^{\circ n} (a_n(s))-\psi ^{\circ n} (b_n(s))\vert}\\
 &\le &
 \vert a_n(s)-b_n(s)\vert \prod _{k=0}^{n-1}\vert \psi '(\psi ^{\circ k} (w))\vert\\
&\le &
\vert a_n(s)-b_n(s)\vert \prod _{k=0}^{n-1}\left( 1+\frac{1}{\lambda _{n-k}^2}\right)\\
&=&
\lambda _n\left( \left(\frac{\lambda _{n}}{\lambda _{n-1}}\right)^s-
\left( \frac{\lambda _{n+1}}{\lambda _{n}}\right) ^s\right) \prod _{k=1}^{n}\left( 1+\frac{1}{\lambda _{k}^2}\right)\\
&\le& \lambda _n\left( \left(\frac{\lambda _{n}}{\lambda_{n-1}}
\right)^s-\left( \frac{\lambda _{n+1}}{\lambda _{n}}\right) ^s\right)
 \prod _{k=1}^{n}\left( 1+\frac{1}{k}\right)\\
&=& (n+1)\lambda _n\left( \left(\frac{\lambda _{n}}{\lambda _{n-1}}\right)^s-
\left( \frac{\lambda _{n+1}}{\lambda _{n}}\right) ^s\right),
\end{eqnarray*}
where we have used $\sqrt{k}\leq\lambda_k$ from Lemma \ref{thm:aps}
part 1.

For $1<y<x$ there exists $y<\xi<x$ such that
\begin{equation}\label{eq:meanvalue}
 x^s-y^s= s(x-y)\xi^{s-1}\le\left\{\begin{array}{ll}
s(x-y) & \mbox{if $0<s\le 1$}\\
sx^{s-1}(x-y) & \mbox{if $1<s$},
\end{array}
\right.
\end{equation}
 and using
$\lambda_n/\lambda_{n-1}\le\lambda_2$ for $n\ge 2$ we get
$$
\left(\frac{\lambda_n}{\lambda_{n-1}}\right)^s-
\left(\frac{\lambda_{n+1}}{\lambda_{n}}\right)^s\le s\max(\lambda_2^{s-1},1)
\left(\frac{\lambda_n}{\lambda_{n-1}}-
\frac{\lambda_{n+1}}{\lambda_{n}}\right).
$$

By (\ref{eq:fes2}) we have
\begin{eqnarray*}
\lefteqn{
 \frac{\lambda _{n}}{\lambda _{n-1}}-
 \frac{\lambda _{n+1}}{\lambda _{n}}}\\
&=&\frac12\left(\sqrt {1+ \frac{4}{\lambda _{n-1}^2}}-
\sqrt {1+ \frac{4}{\lambda _{n}^2}}\,\right)
=\frac{2\displaystyle \left(\frac{1}{\lambda _{n-1}^2}-
\frac{1}{\lambda _{n}^2}\right)}{\sqrt {1+ \frac{4}{\lambda _{n-1}^2}}+
\sqrt {1+ \frac{4}{\lambda _{n}^2}}}
\le \frac{\lambda _n^2-\lambda ^2_{n-1}}{\lambda^2 _n\lambda ^2_{n-1}},
\end{eqnarray*}
so the final estimate is
$$
\vert \psi ^{\circ n} (a_n(s))-\psi ^{\circ n} (b_n(s))\vert
\le \max(\lambda_2^{s-1},1)\frac{s(n+1)}{\lambda _n\lambda ^2_{n-1}}(\lambda _n^2-\lambda ^2_{n-1}),
$$
which tends to zero by part 2, 3 and 4 of Lemma \ref{thm:aps}.

\medskip
\noindent $2^{\circ}:$
{\it Equation \eqref{eq:iterateab} holds for $0< s\le 1$.}

\medskip
This is part of Theorem \ref{thm:main-b}, but for the convenience of
the reader we repeat the main step.

For any convex function $h$ on $]0,\infty[$  we have for $0< s\leq 1$
and $n\ge 2$ 
$$
h(n)-h(n-1)\le \frac{h(n+s)-h(n)}{s}\le h(n+1)-h(n).
$$
By taking $h =\log (1/f)$, which it is convex because $1/f$ is
completely monotonic, we get using $f(n)=\lambda_n$
$$
\log \frac{\lambda _{n-1}}{\lambda _n}\le \frac{1}{s}\log \frac{\lambda_n}{f(n+s)}\le \log \frac{\lambda _{n}}{\lambda _{n+1}},
$$
hence
$$
\lambda _n\left( \frac{\lambda _{n+1}}{\lambda _{n}}\right) ^s
\le f(n+s)\le
\lambda _n \left(\frac{\lambda _{n}}{\lambda _{n-1}}\right) ^s, \quad
0<s\le 1.
$$

Using that $\psi$ is increasing
on $]0,\infty[$, we get by applying $\psi^{\circ n}$  to the previous
 inequality
$$
\psi ^{\circ n} (b_n(s))\le f(s)=\psi^{\circ n}(f(n+s))\le \psi ^{\circ n} (a_n(s)),
$$
and $2^\circ$ follows from $1^\circ$.

\medskip
\noindent$3^{\circ}:$ {\it Equation \eqref{eq:iterateab} holds for $0< s<\infty$.}

\medskip
We will show that if \eqref{eq:iterateab} holds for some $s>0$, then
it also holds for $s+1$. By $2^\circ$ we then get $3^\circ$.

Assume now that \eqref{eq:iterateab} holds for some $s>0$. In particular
$$
\psi^{\circ(n+1)}(a_{n+1}(s))\to f(s)=\psi(f(s+1)).
$$
Using 
\begin{equation}\label{eq:psiinvers}
\varphi(x)=(1/2)(x+\sqrt{x^2+4}),x\in\mathbb R
\end{equation}
as the continuous inverse of $\psi|\;]0,\infty[$
we get
$$
\psi^{\circ  n}(a_{n+1}(s))=\psi^{\circ n}(b_n(s+1))\to\varphi(f(s))=f(s+1),
$$
which is ``half'' of   \eqref{eq:iterateab} for $s+1$, but the other
``half'' comes from $1^\circ$. 

\medskip
Notice that $\lambda_n-\lambda_{n-1}=1/\lambda_n$ and hence
\begin{equation}\label{eq:rho}
\rho_{n,N}=\lambda_n-\lambda_{n-N}=\sum_{k=1}^N\frac{1}{\lambda_{n+1-k}},\quad
n\ge N\ge 1.
 \end{equation}

We denote $D(a,r)=\{z\in\mathbb C\mid |z-a|<r\}$.

\medskip
\noindent$4^{\circ}:$ {\it For $N\in\mathbb N,0<c\le 1,n>N$}
$$
\psi(D(\lambda_n,c\rho_{n,N}))\subseteq D(\lambda_{n-1},c\rho_{n-1,N}).
$$

\medskip
If $|z-\lambda_n|<c\rho_{n,N}$ we get
$$
|\psi(z)-\lambda_{n-1}|=|z-\lambda_n-(\frac{1}{z}-\frac{1}{\lambda_n})|
<c\rho_{n,N}(1+\frac{1}{|z|\lambda_n})
$$
and
$$
|z|=|\lambda_n-(\lambda_n-z)|\ge
\lambda_n-|\lambda_n-z|>\lambda_n-\rho_{n,N}
=\lambda_{n-N},
$$
hence
$$
|\psi(z)-\lambda_{n-1}|<c\rho_{n,N}(1+\frac{1}{\lambda_n\lambda_{n-N}})
=c(\sum_{k=1}^N\frac{1}{\lambda_{n+1-k}}+\frac{\lambda_n-\lambda_{n-N}}
{\lambda_n\lambda_{n-N}})=c\rho_{n-1,N}.
$$

Iterating $n-N$ times using $\rho_{N,N}=\lambda_N-\lambda_0=\lambda_N$
we get

\medskip
\noindent$5^\circ:$ {\it For $1\le N\le n$ and $0<c\le 1$}
$$
\psi^{\circ(n-N)}(D(\lambda_n,c\rho_{n,N}))\subseteq D(\lambda_N,c\lambda_N).
$$

\medskip
\noindent$6^\circ$: {\it For $0<c\le 1, |z|\le cN,N\le n $ we have $b_n(z)\in D(\lambda_n,c\rho_{n,N}).$}

\medskip
For $a>1$ and $z\in\mathbb C, |z|\le 1$ we have the
elementary inequality
\begin{equation}\label{eq:exp}
|a^z-1|\le |z|(a-1).
\end{equation}
Applying this with $a=(\lambda_{n+1}/\lambda_n)^N$ we get
$$
|b_n(z)-\lambda_n|=\lambda_n|(\frac{\lambda_{n+1}}{\lambda_n})^z-1|
\le\lambda_n c\left((\frac{\lambda_{n+1}}{\lambda_n})^N-1\right)
$$
$$
=\lambda_n c(\frac{\lambda_{n+1}}{\lambda_n}
-1)\sum_{k=0}^{N-1}(\frac{\lambda_{n+1}}{\lambda_n})^k
=\frac{c}{\lambda_{n+1}}\sum_{k=0}^{N-1}(\frac{\lambda_{n+1}}{\lambda_n})^k
$$
$$
=c\left(\frac{1}{\lambda_{n+1}}+\sum_{k=1}^{N-1}\frac{\lambda_{n+1}^{k-1}}
{\lambda_n^k}\right)\le
  c\sum_{k=0}^{N-1}\frac{1}{\lambda_{n-k+1}}<
 c\sum_{k=0}^{N-1}\frac{1}{\lambda_{n-k}}=c\rho_{n,N},
$$
where we have used the inequalities
\begin{equation}\label{eq:concave}
\lambda_{n+1}^{k-1}\lambda_{n-k+1}\le\lambda_n^k,\quad k=1,\ldots,N-1,
\end{equation}
which are equivalent to
$$
(k-1)\log f(n+1) + \log f(n-k+1)\le k\log f(n),
$$
but they hold because $\log f$ is concave.

\medskip
Combining $5^\circ$ and $6^\circ$ we get

\medskip
\noindent$7^{\circ}:$ {\it For $0<c\le 1, n\ge N,|z|\le cN$
$$
\psi^{\circ(n-N)}(b_n(z))\in D(\lambda_N,c\lambda_N).
$$
In particular, the sequence $\psi^{\circ(n-N)}(b_n(z)),n\ge N$ of
holomorphic
functions in the disc $D(0,N)$ is bounded on compact subsets of this disc.}

\medskip
\noindent$8^{\circ}:$ {\it For $0< s<\infty, n\ge N$ we have}
$$
\lim_{n\to\infty} \psi^{\circ(n-N)}(b_n(s))=f(s+N).
$$

Since $b_n(s)>\lambda_n$ we know that
$\psi^{\circ(n-N)}(b_n(s))>\lambda_N$. For each $N\ge 1$ we see that
$\psi^{\circ N}|\;]\lambda_{N-1},\infty[\to\mathbb R$ is a
homeomorphism with inverse $\varphi^{\circ N}$, where $\varphi$ is
given by \eqref{eq:psiinvers}. Since
$\psi^{\circ n}(b_n(s))\to f(s)$ we get
$$
\varphi^{\circ N}(\psi^{\circ n}(b_n(s))\to\varphi^{\circ
  N}(f(s))=f(s+N)
$$
i.e.
$$
\lim_{n\to\infty}\psi^{\circ(n-N)}(b_n(s))=f(s+N).
$$ 

\medskip
\noindent$9^{\circ}:$ {\it Let $N\in\mathbb N$. For $z\in D(0,N)$ we have
$$
\lim_{n\to\infty}\psi^{\circ(n-N)}(b_n(z))=f(z+N),
$$
and the convergence is uniform on compact subsets of $D(0,N)$.}

\medskip
By Montel's theorem the sequence
$\psi^{\circ(n-N)}(b_n(z))$ has accumulation points $h$ in
the space $\mathcal H(D(0,N))$ of holomorphic functions on $D(0,N)$.
By $8^\circ$ we know that $h(s)=f(s+N)$ for $0<s<N$.
By the uniqueness theorem for holomorphic functions, all accumulation
points then agree with $f(z+N)\in\mathcal H(D(0,N))$, and the result follows.

\medskip
\noindent$10^{\circ}$: {\it For $z\in\mathbb C$ we have
$$
\lim_{n\to\infty}\psi^{\circ n}(b_n(z))=f(z),
$$
uniformly on compact subsets of $\mathbb C$.}

\medskip
For a compact subset $K\subset\mathbb C$ we choose $N\in\mathbb N$
such that $K\subset D(0,N)$ and know by $9^{\circ}$ that
$\psi^{\circ(n-N)}(b_n(z))$ converges uniformly to $f(z+N)$ for $z\in
K$. We next use that
$\psi^{\circ N}:\mathbb C^*\to\mathbb C^*$ is
continuous, hence uniformly continuous with respect to the chordal
metric on $\mathbb C^*$, and since
$\psi^{\circ N}(f(z+N))=f(z)$, the result follows.

\medskip
\noindent$11^{\circ}:$ {\it For $z\in\mathbb C$ we have
$$
\lim_{n\to\infty}\psi^{\circ n}(a_n(z))=f(z),
$$
uniformly on compact subsets of $\mathbb C$.}

\medskip
In fact,
$$
\psi^{\circ n}(a_{n+1}(z))=\psi^{\circ n}(b_n(z+1))\to f(z+1)
$$
so
$$ 
\psi(\psi^{\circ n}(a_{n+1}(z)))\to \psi(f(z+1))=f(z).
$$

\medskip
{\bf Acknowledgment} The authors want to thank Michael Olesen,
Nyk{\o}bing Katedralskole for
useful comments to an early version of this manuscript.


\begin{thebibliography}{abc}

\bibitem{ak}  Akhiezer,  N.\  I., {\it The classical moment problem}.
Oliver and Boyd, Edinburgh, 1965.

\bibitem{ar}  Artin, E., {\it The Gamma Function}. Holt, Rinehart
  and Winston, New York, 1964.

\bibitem{bear} Beardon,  A. F., {\it Iteration of rational functions}.
 Graduate Texts in Mathematics vol. {\bf 132}.
Springer-Verlag, Berlin-Heidelberg-New York, 1991.

\bibitem{bd}  Berg,  C., Dur\'an, A. J., {\it A transformation
from Hausdorff to Stieltjes moment sequences}. Ark. Mat. {\bf 42} (2004), 239--257.

\bibitem{bd2} Berg, C., Dur\'an, A. J., {\it Some transformations of Hausdorff
moment sequences and Harmonic numbers}. Canad. J.
Math. {\bf 57} (2005), 941--960.

\bibitem{bd3} Berg, C., Dur\'an, A. J., {\it The fixed point for a
 transformation of Hausdorff moment sequences and iteration of a
 rational function}. To appear in Math. Scand.

\bibitem{ha} Hausdorff, F., {\it Momentenprobleme f\"ur ein endliches
 Intervall,} Math. Z. {\bf 16} (1923) 220--248.

\bibitem{W} Widder, D. V., {\it The Laplace Transform}. Princeton
  University Press, Princeton, 1941.

\end{thebibliography}
\end{document}